\documentclass[12pt]{amsart}
\usepackage[usenames,dvipsnames]{xcolor}
\usepackage{tikz}
\usepackage{fancyhdr}
\usepackage{txfonts}
\usepackage{graphicx}
\usepackage{epsfig}
\usepackage{mathrsfs}
\usepackage{amssymb}
\usepackage{latexsym}
\usepackage{amsmath}
\usepackage{amsfonts}
\usepackage{amsbsy}
\usepackage{amscd}
\usepackage{subfigure}
\usepackage{verbatim}
\usepackage{color,tikz}
\usetikzlibrary{calc}

\usepackage{amsthm}
\usepackage{amsfonts}
\usepackage{amsbsy,calc}
\usepackage{graphicx,ifthen,cite}
\usepackage{bm}
\usepackage[normalem]{ulem}

\newcommand{\B}{{\mathcal B}}
\newcommand{\D}{{\mathcal D}}
\newcommand{\CC}{{\mathcal C}}

\newcommand{\R}{{\mathbb R}}
\newcommand{\Z}{{\mathbb Z}}
\newcommand{\C}{{\mathbb C}}

\newcommand{\N}{{\mathbb N}}

\newcommand{\al}{\alpha}

\newcommand{\la}{\lambda}
\newcommand{\La}{\Lambda}

\def\wh{\widehat}

\newtheorem{prop}{Proposition}[section]
\newtheorem{lem}[prop]{Lemma}
\newtheorem{ques}[prop]{Question}
\newtheorem{defi}[prop]{Definition}
\newtheorem{coro}[prop]{Corollary}
\newtheorem{theo}[prop]{Theorem}
\newtheorem{conj}[prop]{Conjecture}
\newtheorem{exam}[prop]{Example}

\newcommand{\ez} {\hspace*{\fill} $\square$}
\newcommand{\eproof}{\hfill$\square$}

\newif\ifdraft
\drafttrue
\ifdraft
\else
\fi
\numberwithin{equation}{section}
\title{}
\author{}
\date{}

\begin{document}
	\baselineskip 16pt
	
	\title[]{{On a distinctive property of Fourier bases associated with $N$- Bernoulli Convolutions}}
	\author[Z.-C. CHi]{Zi-Chao Chi}
	\address{School of Mathematics and Statistics, Key Lab NAA-MOE and Hubei Key Lab-Math. Sci., Central China Normal University, Wuhan 430079, P. R. China}
	\email{skychizichao@163.com}
	\author[X.-G. He]{Xing-Gang He*}
	\address{School of Mathematics and Statistics, Key Lab NAA-MOE and Hubei Key Lab-Math. Sci., Central China Normal University, Wuhan 430079, P. R. China}
	\email{xingganghe@163.com}
	\author[Z.-Y. Wu]{Zhi-Yi Wu}
	\address{School of Mathematics and Information Science, Guangzhou University, Guangzhou, 510006, P. R. China}
	\email{zhiyi\_wu2021@163.com}
	\thanks{*Corresponding author.}
	\thanks{This work was supported by the National Natural Science Foundation of China 12371087 and 12301105}
\begin{abstract}
A distinctive problem of harmonic analysis on $\R$ with respect to a Borel probability measure $\mu$ is identifying all $t\in\R$ such that both
\[\left\{e^{-2\pi i\lambda x}: \lambda\in\Lambda\right\}\quad\text{and}\quad \left\{e^{-2\pi i\lambda x}: \lambda\in t\Lambda\right\}\]
form orthonormal bases of the space $L^2(\mu)$. Currently, this phenomenon has been observed only in certain singular measures. It is deeply connected to the convergence of Mock Fourier series with respect to the aforementioned bases. In this paper, we apply classical number theory to solve the general conjecture and basic problems in this field within the setting of $N$-Bernoulli convolutions, which extend almost all known results and give some new ones.
\end{abstract}

\maketitle
\tableofcontents

\section{Introduction}
Let $\mu$ be a Borel probability measure with compact support in $\R^d$, and let $\Lambda \subset \R^d$ be a countable set. The pair $(\mu, \Lambda)$ is called a {\it spectral pair} if the set
\begin{eqnarray*}
	E_\La=\left\{e^{-2\pi i \langle\la,  x\rangle}: \la\in \La\right\}
\end{eqnarray*}
forms an orthonormal basis of $L^2(\mu)$. In this case, $\mu$ is referred to as a {\it spectral measure}, and $\Lambda$ as a {\it spectrum} of $\mu$. He, Lai, and Lau \cite{HLL2013} proved that a spectral measure $\mu$ in $\R^d$ is always of pure type. That is, $\mu$ is either a finite sum of Dirac measures with equal probability weights or is singular or absolutely continuous with respect to (w.r.t.) the Lebesgue measure.

It is well established that an absolutely continuous spectral measure must correspond to the normalized Lebesgue measure on its support \cite{DL14}. The spectrality of such a measure, $\mu = \frac{1}{\mathcal{L}(\Omega)} \mathcal{L}_{\Omega}$, where $\mathcal{L}_{\Omega}$ denotes the Lebesgue measure restricted to a set $\Omega$, depends on whether $\Omega$ satisfies the translation tile property in $\R^d$. This connection was initially conjectured by Fuglede \cite{Fu74} and has since been investigated in numerous subsequent studies.

A classical example of a spectral measure is the Lebesgue measure $\mathcal{L}_{[0, 1]}$, which has the unique spectrum $\Z$ up to translations. In contrast, a singular spectral measure may have infinitely many spectra, also up to translations. For instance, let  $\mu_{4, \{0, 2\}}$ be the middle forth Cantor measure ($\mu_{4, \{0, 2\}}$ is defined in the same way as in \eqref{eqmeasure} below) and
$$\Lambda(4, \{0, 1\})=\left\{\sum_{k=0}^n 4^k c_k: \text{all $c_k\in \{0, 1\}$ and $n\ge 0$}\right\}.$$
Then, $\left( \mu_{4, \{0, 2\}}, \Lambda(4, \{0, 1\})\right) $ is a spectral pair (while the middle third Cantor measure $\mu_{3, \{0, 2\}}$ is not a spectral one), as proved by Jorgenson and Pedersen \cite{JP1998}. Moreover, there are infinitely many
$$t\in\N^\ast=\{1, 2, 3, \ldots\}$$
such that  $(\mu_{4, \{0, 2\}}, t\Lambda(4, \{0, 1\})$ is also a spectral pair  \cite{DH16}.

As applications of the above results, Strichartz \cite{Str2006} showed that the mock Fourier series of a continuous function $f$ w.r.t. the spectral pair $(\mu_{4, \{0, 2\}}, \Lambda(4, \{0, 1\}))$ converges uniformly to $f$ on the support of $\mu_{4,\{0,2\}}$. On the other hand, Dutkay et al.\cite{DHS2014} found that there exists a continuous function whose Fourier series w.r.t. the spectral pair $(\mu_{4, \{0, 2\}}, t\Lambda(4, \{0, 1\})$, $t=17$ or 23,  diverges at $0$. Those surprising phenomenons motivate the theme of this paper. We begin by defining (following \cite{DH16})
\begin{defi}
	The $t\in\R$ is called a complete (or an incomplete) number of the spectral pair $(\mu, \Lambda)$ if the pair $(\mu, t\Lambda)$ is (or is not) a spectral pair, respectively.
\end{defi}
It is worth noting that $t$ is referred to as {\it eigenvalue} or {\it scaling number} in some other works. The phenomenon of complete number of a spectral pair was first discovered  by Strichartz\cite[Example 2.9(b)]{Str2000}  and independently by {\L}aba and Wang  \cite[Example 3.2]{LW2002}.

In general, there is a general conjecture as follows:
\begin{conj}\label{conj1}
	Let  $(\mu, \Lambda)$ be a singular spectral pair in $\R$. Then there are infinitely many complete numbers of this pair.
\end{conj}
{If this conjecture holds true, then a natural question arises:}
\begin{ques}\label{q1}
	Let $p$ be a prime. Is $p^n$ complete of a given singular spectral pair $(\mu, \Lambda)$  for each $n\ge 1$?
\end{ques}

Recent papers \cite{DH16} and \cite{DK18} have investigated the Cantor measure $\mu_{2q, \{0, q\}}$ for $q=2$ and $q\ge 2$, respectively. These works provide important contributions to the complete number problem. In particular, Dutkay and Haussermann \cite{DH16} discovered that the complete number problem is closely related to classic  number theory. In order to prepare  the study of Fourier analysis on Cantor measures, which is an extension of trigonometric function theory, in this paper, we will extend their main results to a generalized case by a different method and get some more deeper understanding on the topic.

It is difficult to completely settle Conjecture \ref{conj1} and Question \ref{q1} in general case. In this paper, we study the measure defined by the infinite convolution of finite sums of Dirac measures with equal weights, that is,
\begin{eqnarray}\label{eq1.1}
	\mu_{\rho, \D}=\delta_{\rho^{-1}\D}\ast\delta_{\rho^{-2}\D}\ast\cdots\ast\delta_{\rho^{-n}\D}\ast\cdots,
\end{eqnarray}
where $\rho>1$, $\D=\{0,1,\ldots,N-1\}$ for some $N\in \N^*$ and the measure $\delta_E=\frac 1{\#E}\sum_{e\in E}\delta_e$ for any finite set $E$ is a finite sum of Dirac measures with equal weights. The infinite convolutions is understood to converge in the weak sense. Such a measure is a natural generalization of the Bernoulli convolution (in which $\D=\{0,1\}$), and we refer to it as the $N$-Bernoulli convolution. The close relation to number theory makes the study of the case of a general integer $N$ particularly challenging.

We now recall the definition of {\it self-similar measures} (see, e.g., \cite{F1990}).  For a real number $\rho>1$ and finite set $\B\subset\R$, the associated self-similar measure $\mu_{\rho, \B}$ is the unique Borel probability measure which satisfies the following equation:
\begin{eqnarray}\label{eqmeasure}
	\mu_{\rho, \B}(\cdot)=\frac 1{\#\B}\sum_{b\in \B}\mu_{\rho, \B}\circ \phi_b^{-1}(\cdot),
\end{eqnarray}
where $\{\phi_b(x)=\rho^{-1}(x+b):b\in\B\}$ forms an {\it iterated function system (IFS)}. The support of $\mu_{\rho, \B}$ is the set
\begin{eqnarray*}
	T(\rho, \B)=\left\{\sum_{k=1}^\infty \rho^{-k}b_k: \text{all $b_k\in \B$}\right\},
\end{eqnarray*}
which is a fractal set.

If $\B=\D$, then $\mu_{\rho, \B}$ is the $N$-Bernoulli convolution $\mu_{\rho, \D}$. Notably, the cases $\mu_{N, \D}=\mathcal{L}_{[0, 1]}$ and $\mu_{\rho, \{0\}}=\delta_0$, respectively.
Furthermore, it is known that the $\mu_{\rho, \D}$ is singular continuous w.r.t. Lebesgue measure if $\rho>\#\B$, and it is a spectral measure if and only if $\rho\in\N^*$ and $N\mid \rho$, as proved in \cite{JP1998, HL2008, Dai2012, DHL2013, Deng2014, DHL2014}. This model forms the foundation of numerous analytical theories, as seen in the works \cite{Dai2016, DHL2013,FH2017,HTW2018,WW24} and references therein.

In this paper, we adopt a new notation $b$ in replace of $\rho$, where $b=qN$ with $q\ge 2$ being an integer. We define
\begin{eqnarray}\label{spectrum}
	\Lambda(b, \CC)=\bigcup_{k=1}^\infty \Lambda_k(b, \CC)\quad\text{with}\quad \Lambda_k(b, \CC)=\CC\oplus b\CC\oplus\cdots\oplus b^{k-1}\CC,
\end{eqnarray}
where $$\CC=\left(-\frac{N}{2},\frac{N}{2}\right]\cap\Z.$$
It follows from \cite[Page 231]{Str1998} that the pairs $(\mu_{b, \D}, q\Lambda(b, \CC))$ and $(\mu_{b, q\D}, \Lambda(b, \CC))$ form spectral pairs (see Lemma \ref{lemspectrum} in Section 2). In this study, we focus on the spectral pair $(\mu_{b, q\D}, \Lambda(b, \CC))$, noting that the complete number problem of the spectral pair $(\mu_{b, q\D}, \Lambda(b, \CC))$ is equivalent to that of $(\mu_{b, \D}, q\Lambda(b, \CC))$ (see Lemma \ref{lem2.4}).

Throughout the remainder of this paper, we simply refer to a complete number of the spectral pair $(\mu_{b, q\D}, \Lambda(b, \CC))$ as a complete number.
The first result of this paper concerns the density of complete numbers. Noticing that the complete number status of $t$ and $-t$ is identical (can be obtained directly from Theorem \ref{theo2.1}), it suffices to restrict our focus on the positive number. It is well known that $\pi(x)\sim{x/\log x}$ ($x\rightarrow+\infty$), where $\pi(x)$ is the prime-counting function. And there are primes which are incomplete numbers (see Example \ref{exam2.12} and \ref{exam2.13}). Using the result from Erd\"{o}s and Murty \cite{EM1999}, we derive the following bound on the density of complete numbers:
\begin{theo}\label{thmdeig}
	Let $\Gamma$ be the set of all positive complete numbers. Then
	$$\limsup_{x\to\infty}\frac{\#(\Gamma\cap[0,x])}{x}\le c\cdot\frac{\varphi(N)}{N},$$
	where $\varphi$ is Euler function and $0<c<1$ is a constant depending on $b$. Moreover, if $q> N$, then
	$$\liminf_{x\to\infty}\frac{\#(\Gamma\cap[0,x])}{\pi(x)}\ge 1.$$
\end{theo}
Proposition \ref{prop2.3} in Section 2 says that a complete number $t$ must be integer with $\gcd(t, N)=1$. The above theorem shows that not all integers $t$ coprime with $N$ are complete numbers, and that infinitely many complete numbers exist provided $q >N$. Moreover, it can be shown that infinitely many complete numbers exist even without this assumption (Theorem \ref{theoinfinite}).

The following theorem allow us to focus only on the positive integer  $t$ that $\gcd(t,b)=1$.
\begin{prop}\label{theoprimitive}
	Any $t\in\N^*$ that $\gcd(t,N)=1$ is a complete number if and only if $\frac{t}{\gcd(t,b)}$ is a complete number.
\end{prop}
Using this fact, we now introduce the following definition to support the further study:
\begin{defi}
	Suppose that $t\in\N^\ast$ and $\gcd(t,N)=1$.
	\begin{enumerate}
		\item A complete number $t$ is called primitive if $\gcd(t, b)=1$.
		\item An incomplete number $t$ is called primitive if any proper divisor $d$ of $t$, $d$ is a complete number.
	\end{enumerate}
\end{defi}
If $t$ is a primitive incomplete number, then $\gcd(t,b) = 1$ (see Theorem \ref{theo1.10}). The primitive number $t$ (whether it is a complete number or an incomplete number) plays a crucial role in understanding the spectrality of the pair $(\mu_{b, q\D}, t\Lambda(b, \CC))$.
With loss of generality we focus particularly on primitive numbers, as the condition $\gcd(t,b) = 1$ allows the applications of number theory techniques.

\begin{theo}\label{theoinfinite}
	There are infinitely many primitive complete numbers and primitive incomplete numbers.
\end{theo}

Even for the spectral pair ($\mu_{b,{q\D}},\Lambda(b,\CC))$, the answer to Question \ref{q1} may not be affirmative. Dutkay and Kraus provided an example \cite[Remark 2.14]{DK18}, showing that the prime $55987$ is an incomplete number of the spectral pair $(\mu_{6, \{0, 3\}}, \Lambda(6, \{0, 1\}))$. In this paper, we provide another example, detailed in Example \ref{exam2.12}. The underlying reason for this are unclear. Fortunately, we can solve Question \ref{q1} under certain conditions, as proved by the following theorem in Section \ref{sectwo}:
\begin{theo}\label{theo1.7}
	Let $r\ge 2$ be an integer and let $p\ne 2^r-1$ be an odd prime. If $rO_2(p)$ is even or $r\nmid(p-1)$, then $p^n$ is a complete number of the spectral pair $(\mu_{2^r, \{0, 2^{r-1}\}}, \Lambda(2^r, \{0, 1\}))$ for each $n\ge 1$.
\end{theo}
We remark that the case $r = 2$ in Theorem \ref{theo1.7} was previously proved by Dutkay and Haussermann in \cite{DH16} using a different method. Additionally, the number $p = 2^r - 1$ is an incomplete number (see Corollary \ref{b-1}).

The following theorem generalizes the results from \cite{DH16, DK18} and provides a sufficient condition for composite numbers to become complete numbers:
\begin{theo}
	Let $p_1,\ldots,p_m\ge 2N+1$ be distinct primes with $\gcd(p_1\cdots p_m,b)=1$. If $p_1^{i_1}\cdots p_m^{i_m}$ is a complete number for $i_1,\ldots,i_m$ sufficiently large, then  $p_1^{k_1}\cdots p_m^{k_m}$ is a complete number for any $k_1,\ldots,k_m\ge 0$.
\end{theo}

The precise meaning of ``sufficiently large" can be found in Theorem \ref{theodut}.

We organize this paper as follows. In Section \ref{secpre}, we present some preliminary results, simple facts, and the proof of Theorem \ref{theo2.5}. Section \ref{secdist} is dedicated to the density of the complete numbers. Section \ref{secprconj} focuses on the study of primitive numbers, including both complete numbers and incomplete numbers, and provides the proof of Conjecture \ref{conj1} in the case of this paper. In Section \ref{secsuffi}, we introduce sufficient conditions for composite numbers to be complete numbers. Finally, in Section \ref{sectwo}, we give a positive answer to Question \ref{q1} in the case where $N = 2$, $b = 2^r$ for $r > 2$, which generalizes the result of Dutkay et al. in \cite{DH16}.

\bigskip

\section{Preliminaries and some simple facts}\label{secpre}

\bigskip

Let $\mu$ be a Borel probability measure with compact support in $\R$. The Fourier transform of $\mu$ is defined as
\[\wh{\mu}(\xi)=\int e^{-2\pi i \xi x} d\mu(x).\]
The zero set of a function $f$ is denoted by $\mathcal{Z}(f)=\{x: f(x)=0\}$. Let $\Lambda\subset \R$ be a countable set. The set $\Lambda$ is called an {\it orthogonal set} of $\mu$ if the family $E_\Lambda$ forms an orthogonal set in $L^2(\mu)$. This is equivalent to the condition that
\[\Lambda-\Lambda\subset \mathcal{Z}(\wh{\mu})\cup\{0\}.\]
To establish the orthogonality or spectrality of the pair $(\mu, \Lambda)$, we will rely on the following criterion, which was proved by Jorgensen and Pedersen in \cite[Lemma 3.3]{JP1998} and is derived from Parseval's identity.
\begin{theo}\label{theo2.1}
	Let $\mu$ be a Borel probability measure with compact support in $\R$ and denote
	$$Q_\Lambda(\xi)=\sum_{\lambda\in\Lambda}|\wh{\mu}(\xi+\lambda)|^2,\quad\forall \xi\in\R.$$
	Then the following statements hold:
	\begin{enumerate}
		\item [(i)]The set $\Lambda$ is an orthogonal set of $\mu$ if and only if $Q_\Lambda(\xi)\le 1$ for all $\xi\in\R$;
		\item [(ii)]The set $\Lambda$ is a spectrum of $\mu$ if and only if $Q_\Lambda(\xi)\equiv 1$ for $\xi\in\R$.
	\end{enumerate}
	In the above two cases, the function $Q_{\Lambda}$ has an entire analytic extension to $\C$.
\end{theo}

Let $\mu_{b, q\D}$ be defined by \eqref{eqmeasure}. Then, its Fourier transform is given by
\[\wh{\mu}_{b, q\D}(\xi)=\prod_{k=1}^\infty M_\D(b^{-k}q\xi),\]
where $M_\D(\xi)$ is the {\it mask function} of $\D$, defined as
\[M_\D(\xi)=\frac 1N\sum_{d\in \D}e^{-2\pi i d\xi}.\]
It follows that
\begin{eqnarray*}\label{eq2.1}
	\mathcal{Z}(M_\D)=N^{-1}(\Z\setminus N\Z)\quad\text{and}\quad \mathcal{Z}(\wh{\mu}_{b, q\D})=\bigcup_{k=1}^\infty b^kq^{-1}\mathcal{Z}(M_\D)=\bigcup_{k=0}^\infty b^k(\Z\setminus N\Z).
\end{eqnarray*}
Let $\Lambda(b,\CC)$ be defined as in \eqref{spectrum}. Then we have the following lemmas:
\begin{lem}\label{lem2.4}
	$(\mu_{b, \D}, tq\Lambda(b, \CC))$ is a spectral pair if and only if $(\mu_{b, q\D}, t\Lambda(b, \CC))$ is a spectral pair.
\end{lem}
\proof According to Theorem \ref{theo2.1}, the assertion follows from
\begin{eqnarray*}
	1&\equiv& \sum_{\lambda\in\La(b, \CC)}|\wh{\mu}_{b, q\D}(\xi+t\la)|^2=\sum_{\lambda\in\La(b, \CC)}\prod_{k=1}^\infty|M_{q\D}(b^{-k}(\xi+t\la))|^2\\
	&=&\sum_{\lambda\in\La(b, \CC)}\prod_{k=1}^\infty|M_{\D}(b^{-k}(q\xi+tq\la))|^2= \sum_{\lambda\in\La(b, \CC)}|\wh{\mu}_{b, \D}(q\xi+tq\la)|^2.
\end{eqnarray*}
\ez

\begin{lem}\label{lemspectrum}
	Both $(\mu_{b, q\D}, \Lambda(b,\CC))$ and  $(\mu_{b, \D}, q\Lambda(b,\CC))$ are spectral pairs in $\R$.
\end{lem}
\proof  By Lemma \ref{lem2.4}, we only need to prove the pair $(\mu_{b, q\D}, \Lambda(b, \CC))$ is spectral. We prove it by the theorem on page 231 of \cite{Str1998}, which says that the pair $(\mu_{b, q\D}, \Lambda(b, \CC))$ is a spectral pair if $\mathcal{Z}(M_{b^{-1}q\D})\cap T(b, \CC)=\emptyset$. Clearly $T(b, \CC)\subset \left[-\frac{N}{2}, \frac N{2}\right]$ and $\mathcal{Z}(M_{b^{-1}q\D})=\Z\setminus N\Z$, then the Strichartz's criterion holds and the assertion follows.\ez

The following proposition indicates that we only need to focus on integers that are coprime with $m$.
\begin{prop}\label{prop2.3}
	If $t\in\R$ is a complete number, then $t\in\Z$ and $\gcd(t, N)=1$.
\end{prop}
\begin{proof}
	If $t$ is a complete number, then $t\in t\Lambda(b,\CC)\subset\mathcal{Z}(\widehat{\mu}_{b,{q\D}})\cup\{0\}\subset\Z$. We claim that $t\CC$ is a spectrum of $\delta_{b^{-1}q\D}$. Indeed, by Theorem \ref{theo2.1}, for $\xi\in\R$,
	\begin{equation*}
		\begin{split}
			1&\equiv Q_{t\Lambda(b,\CC)}(\xi)=\sum_{\lambda\in{\Lambda(b,\CC)}}\left\vert\widehat{\mu}_{b,q\D}(\xi+t\lambda)\right\vert^2=\sum_{c\in \CC,\lambda\in{\Lambda(b,\CC)}}\left\vert\widehat{\mu}_{b,q\D}(\xi+tc+tb\lambda)\right\vert^2\\
			&=\sum_{c\in \CC}\sum_{\lambda\in{\Lambda(b,\CC)}}\left\vert M_{q\D}(b^{-1}\xi+t\lambda+tb^{-1}c)\right\vert^2 \prod_{k=2}^{\infty}\left\vert M_{q\D}(b^{-k}(\xi+tc+tb\lambda))\right\vert^2\\
			&=\sum_{c\in \CC}\left\vert M_{q\D}(b^{-1}(\xi+tc))\right\vert^2\sum_{\lambda\in{\Lambda(b,\CC)}}\prod_{k=1}^{\infty}\left\vert M_{q\D}(b^{-k}(b^{-1}\xi+tb^{-1}c+t\lambda))\right\vert^2\\
			&=\sum_{c\in \CC}\left\vert\widehat{\delta}_{b^{-1}{q\D}}(\xi+tc)\right\vert^2.
		\end{split}
	\end{equation*}
	This implies that $t\CC$ is a spectrum of $\delta_{b^{-1}{q\D}}$.
	
	Suppose that $\gcd(t,N)=r>1$. We denote $t=rt'$ and $N=rN'$ for some nonzero $t',N'\in\Z$. Note that $N'\in \CC$ and
	$$tN'=rt'N'=Nt'\notin\Z\setminus N\Z=\mathcal{Z}(\widehat{\delta}_{b^{-1}{q\D}}).$$
	This contradicts $t\CC\setminus\{0\}\subset\mathcal{Z}(\widehat{\delta}_{b^{-1}{q\D}})$ since $t\CC$ is a spectrum of $\delta_{b^{-1}{q\D}}$.
\end{proof}

Based on the above facts, we provide a criterion to determine whether an integer $t$, which is coprime with $N$, is a complete number:
\begin{theo}\label{theo2.5}
	Let $t\in\N^\ast$ with $\gcd(t, N)=1$. Then the following statements are equivalent:
	\begin{enumerate}
		\item[(i)]$t$ is an incomplete number.
		\item[(ii)] There exists  a non-zero integer periodic sequence $\{x_k\}_{k=0}^\infty $ w.r.t. the pair $(b, t\CC)$, i.e., there exist $\{c_k\}_{k=0}^\infty$ in $\CC$ such that
		$x_{k+1}=\frac{x_k+tc_k}b\in\Z\setminus\{0\}$ and $x_{n+k}=x_k$ for $k\in\N$ and some $n\ge 1$.
		\item[(iii)] There exist $c_0, c_1, \ldots, c_{n-1}\in \CC$, not all equal to $0$, such that
		\begin{eqnarray*}
			t\frac{c_0+bc_1+\cdots+b^{n-1}c_{n-1}}{b^{n}-1}\in T(b, t\CC)\cap\Z\setminus\{0\}.
		\end{eqnarray*}
		\item[(iv)] $T(b, t\CC)\cap\Z\ne\{0\}$.
		\item[(v)] $qT(b, t\CC)\cap\Z\ne\{0\}$.
	\end{enumerate}
\end{theo}

\begin{proof}
(i)$\Longrightarrow$(ii). According to \cite[Theorem 1.3]{LW2002}, the pair $(\mu_{b, \D}, tq\Lambda(b, \CC))$ is not a spectral pair if and only if there exists a non-zero integer periodic sequence $\{x_k\}_{k=0}^{\infty}$ w.r.t. $(b, tq\CC)$, that is, there exists $\{c_k\}_{k=0}^\infty$ in $\CC$ such that
\[x_{k+1}=\frac{x_k+tqc_k}b\in\Z\setminus\{0\}, \qquad \text{for $k\in\N$}.\]
Recall that $b=qN$. Thus $q$ divides $x_k$ for $k\ge1$. Let $y_k=q^{-1}x_k$. Then $\{y_k\}_{k=0}^{\infty}$ is a non-zero integer periodic sequence w.r.t. $(b, t\CC)$.

(ii)$\Longrightarrow$(iii). For $k\in\N$, we have that
\begin{eqnarray*}
	y_{k+1}=\frac{y_k+tc_k}b=\frac{y_{k-1}+tc_{k-1}+tbc_k}{b^2}=\cdots=\frac{y_0+tc_0+tbc_1+\cdots+tb^kc_k}{b^{k+1}}.
\end{eqnarray*}
Then since $y_n=y_0$, where $n$ is the minimal positive period,
\begin{eqnarray}\label{eq2.2}
	y_0=t\frac{c_0+bc_1+\cdots+b^{n-1}c_{n-1}}{b^n-1}\in T(b, t\CC)\cap\Z.
\end{eqnarray}
If $y_0=0$, then the above implies that $c_k=0$ for $0\le k\le n-1$, which yields all $x_k=0$. Then (iii) follows.

(iii)$\Longrightarrow$(iv)$\Longrightarrow$(v) is trivial.

(v)$\Longrightarrow$(i).  Notice that
$$qT(b,t\CC)=\bigcup_{c\in \CC}\frac{qT(b,t\CC)+tqc}{b}.$$
Let $x_0\in qT(b,t\CC)\cap\Z\setminus\{0\}$. Then $x_0=\frac{x_{-1}+tqc_{-1}}{b}$ for some $x_{-1}\in qT(b,t\CC)$ and $c_{-1}\in \CC$.
It follows that $x_{-1}=bx_0-tqc_{-1}\in\Z\setminus\{0\}$. Indeed, if $x_{-1}=0$, then $\frac{tc_{-1}}{N}=x_0\in\Z\setminus\{0\}$. Since $\gcd(t,N)=1$ and $c_{-1}\neq 0$, it follows that $c_{-1}$ must be divisible by $N$. This contradicts the fact that $c_{-1}\in\CC=(-\frac{N}{2}, \frac{N}{2}]$.
By induction, we have $x_{-i}\in qT(b,t\CC)\cap\Z\setminus\{0\}$ for any $i\ge 1$. Since $qT(b,t\CC)\cap\Z\setminus\{0\}$ is finite, then there exist $j,k\in\N^\ast$ with $j>k$ such that $x_{-j}=x_{-k}$. This implies that there exists a non-zero integer periodic sequence w.r.t. $(b, tq\CC)$. Then (i) follows from \cite[Theorem 1.3]{LW2002} and Lemma \ref{lem2.4}.
\end{proof}

Based on the above theorem, we can derive several corollaries:
\begin{coro}\label{b-1}
	$b-1$ is an incomplete number. If $N\ge 3$, then $b+1$ is also an incomplete number.
\end{coro}
\begin{proof}
	It is easy to check that $\{1\}_{k=0}^{\infty}$ is a non-zero integer periodic sequence w.r.t. $\left(b,(b-1)\CC\right)$ and $\{(-1)^k\}_{k=0}^{\infty}$ is a non-zero integer periodic sequence w.r.t. $\left( b,(b+1)\CC\right)$ if $N\ge 3$. Hence the assertion follows from (ii) of Theorem \ref{theo2.5}.
\end{proof}


\begin{coro}\label{coro2.10}
	Let $t\in\N^\ast$ with $\gcd(t,N)=1$. If $1\le t< \frac{2(b-1)}N$, then $t$ is a complete number.
\end{coro}
\proof
Note that
\begin{equation*}
	T(b, t\CC)=t\sum_{k=1}^\infty b^{-k}\CC\subset\left( -\frac{tN}{2(b-1)}, \frac{tN}{2(b-1)}\right].
\end{equation*}
For $1\le t<\frac{2(b-1)}{N}$, $T(b,t\CC)\subset (-1,1)$. Then $T(b,t\CC)\cap\Z=\{0\}$ and $t$ is a complete number by (iv) of Theorem \ref{theo2.5}.\ez

\begin{coro}\label{coronn}
	If $t\in\N^\ast$ is an incomplete number, then so is $kt$ for each $k\in\N$.
\end{coro}
\proof If $\gcd(kt, N)>1$, then the assertion follows form Proposition \ref{prop2.3}. If $\gcd(kt, N)=1$, then the assertion follows from (iii) of Theorem \ref{theo2.5}.
\ez

\begin{coro}\label{coro2.11}
	Let $t\in\N^\ast$ with $\gcd(t,N)=1$ be an incomplete number and $\{x_k\}_{k=0}^{\infty}$ be a non-zero integer periodic sequence w.r.t. $(b,t\CC)$. Suppose that $n$ is the minimal positive period. Then
	$$x_k=t\frac{c_{k-n}+bc_{k-n+1}+\cdots+b^{n-1}c_{k-1}}{b^{n}-1}\in T(b, t\CC)\cap\Z\setminus\{0\},\quad k\ge n$$
	for some $\{c_k\}_{k=0}^{\infty}$ in $\CC$ with the same minimal period $n$.
\end{coro}
\begin{proof}
	The proof of the assertion is similar to the proof of (ii)$\Longrightarrow$ (iii) of Theorem \ref{theo2.5} and we omit it.
\end{proof}
\begin{coro}\label{coro2.12}
	Let $t\in\N^\ast$ with $\gcd(t,N)=1$ be an incomplete number. If $\{x_k\}_{k=0}^\infty$ and $\{y_k\}_{k=0}^\infty$ be two distinct non-zero integer periodic sequence w.r.t. $(b,t\CC)$, then $x_i \neq y_j$ for any $i,j\in\N$.
\end{coro}
\begin{proof}
	Let $n$ and $l$ be minimal positive period of $\{x_k\}_{k=0}^\infty$ and $\{y_k\}_{k=0}^\infty$, respectively. Without loss of generality, we suppose $x_0=y_0$. Similar as the proof of (ii)$\Longrightarrow$ (iii) of Theorem \ref{theo2.5}, there exist sequences $\{c_k\}_{k=0}^\infty$ and $\{d_k\}_{k=0}^\infty$ in $\CC$ such that
	\begin{eqnarray*}
		x_0=t\frac{c_0+bc_1+\cdots+b^{n-1}c_{n-1}+\cdots+b^{nl-1}c_{nl-1}}{b^{nl}-1}
	\end{eqnarray*}
	and
	\begin{eqnarray*}
		y_0=t\frac{d_0+bd_1+\cdots+b^{n-1}d_{n-1}+\cdots+b^{nl-1}d_{nl-1}}{b^{nl}-1}.
	\end{eqnarray*}
	$x_0=y_0$ forces $c_i=d_i$ for $0\le i\le nl-1$. Then $x_i=y_i$ for all $i\in\N$, which leads a contradiction.
\end{proof}

\begin{coro}\label{coro2.9}
	Let $t\in\N^\ast$ with $\gcd(t,N)=1$ be a primitive incomplete number. Then $\gcd(t, b)=1$ and $\gcd(t, x_k)=1$ for all $x_k$ in a non-zero integer periodic sequence  $\{x_k\}_{k=0}^\infty$ w.r.t. $(b, t\CC)$.
\end{coro}
\proof Suppose $\gcd(t, b)=r>1$. Then $\gcd(r, b^n-1)=1$ for each $n\ge 1$. Write $t=t'r$. According to (iii) of Theorem \ref{theo2.5}, there exist $c_0, c_1, \ldots, c_{n-1}\in \CC$ such that
\begin{eqnarray*}
	t'\frac{c_0+bc_1+\cdots+b^{n-1}c_{n-1}}{b^{n}-1}\in T(b,t'\CC)\cap\Z\setminus\{0\},
\end{eqnarray*}
and therefore $t'$ is an incomplete number, which contradicts the definition of primitive incomplete number $t$. Hence $\gcd(t, b)=1$.

Let $\{x_k\}_{k=0}^{\infty}$ be a non-zero integer periodic sequence w.r.t. $(b, t\CC)$ and $n$ be its minimal positive period. Without loss of generality, suppose that $\gcd(t,x_0)=r>1$. Write $x_0=rx_0'$ and $t=rt'$. By Corollary \ref{coro2.11} we have
\[x_0=x_n=t\frac{c_0+bc_1+\cdots+b^{n-1}c_{n-1}}{b^{n}-1}\in T(b, t\CC)\cap\Z\setminus\{0\}.\]
Similarly, we have that $t'$ is an incomplete number, which is impossible and the desired results follows.\ez

\begin{theo}\label{theo2.11}
	
	Suppose $t$ is a complete number. If $r$ is a divisor of $b$ with $\gcd(r,N)=1$, then $tr^k$ is a complete number for any $k\in\N$.
\end{theo}
\begin{proof}
	$\gcd(tr^k,N)=1$ follows from Proposition \ref{prop2.3} and $\gcd(r,N)=1$. Suppose $tr^k$ is an incomplete number.  According to Theorem \ref{theo2.5}, there exist $c_0, c_1, \ldots, c_{n-1}\in \CC$, not all equal to $0$, such that
	\begin{eqnarray*}
		tr^k\frac{c_0+bc_1+\cdots+b^{n-1}c_{n-1}}{b^{n}-1}\in T(b, tr^k\CC)\cap\Z\setminus\{0\}.
	\end{eqnarray*}
	Noticing that $\gcd(r,b^n-1)=1$, we have
	\begin{eqnarray*}
		t\frac{c_0+bc_1+\cdots+b^{n-1}c_{n-1}}{b^{n}-1}\in T(b, t\CC)\cap\Z\setminus\{0\}.
	\end{eqnarray*}
	This contradicts that $t$ is a complete number.
\end{proof}
\begin{theo}
	Any positive integer $t$ that $\gcd(t,N)=1$ is a complete number if and only if $\frac{t}{\gcd(t,b)}$ is a complete number.
\end{theo}
\begin{proof}
	It can be obtained by Corollary \ref{coronn} and Theorem \ref{theo2.11}.
\end{proof}

\section{The density of complete numbers}\label{secdist}

 To simplify the analysis, we adopt the following stipulation: throughout this paper, for any $a\in\Z$,
\begin{equation}
	a\pmod t\in \left( -\frac t2, \frac t2\right].
\end{equation}
The following definition is a fundamental concept in number theory, which will be frequently used in this paper:
\begin{defi}
	Let $t\in\N^\ast$ with $\gcd(t,b)=1$. Then the order of $b$ modulo $t$ is defined by
	$$O_b(t)=\min\{\alpha\in\Z:b^{\alpha}\equiv 1\pmod t\}.$$
\end{defi}

\begin{lem}\label{lem4.9}
	Let $t\in\N$ such that $\gcd(t, N)=1$. Then
	\begin{eqnarray*}
		\#\left(T(b,t\CC)\cap\Z\setminus\{0\}\right) <\min_{k\in\N}\left\{\frac{t(N-1)}{q^k(b-1)}+N^k\right\}.
	\end{eqnarray*}
\end{lem}
\proof Recall $\CC=(-\frac N2, \frac N2]\cap\Z$. Then there exists $\alpha\in\Z$ such that $\alpha+\CC=\{0, 1, \ldots, N-1\}$. By the definition of $T(b, t\CC)$ it is easy to obtain that
\[T(b, t\CC)+\frac{t\alpha}{b-1}=T(b, t(\alpha+\CC))\subset\left[0, \frac{t(N-1)}{b-1}\right].\]
Denote $\phi_i(x)=\frac{x+ti}b$, $0\le i\le N-1$. Then for any $I=i_1i_2\cdots i_k\in\{0, 1,\ldots, N-1\}^k$ we have
\[\phi_I(x)=\phi_{i_1}\circ\phi_{i_2}\circ\cdots\circ\phi_{i_k}(x)=\frac x{b^k}+tc_I,\]
where $c_I=\frac {i_1}b+\cdots+\frac{i_k}{b^k}$. According to that
\begin{eqnarray*}
	T(b, t(\alpha+\CC))&=&\bigcup_{I\in\{0, 1, \ldots, N-1\}^k}\phi_I(T(b, t(\alpha+\CC)))\\
	&\subset&\bigcup_{I\in\{0, 1, \ldots, N-1\}^k}\left(tc_I+\left[0, \frac{t(N-1)}{b^k(b-1)}\right]\right),
\end{eqnarray*}
where the union on the right is disjoint, we claim that
\[\frac{t(m-1)}{b^k(b-1)}\notin\Z.\]
In fact, if $\frac{t(N-1)}{b^k(b-1)}\in\Z$, then $\frac{t(N-1)}{N}\in\Z$, which is impossible
by $\gcd(t, N)=1$. Then the interval $tc_I-\frac{t\alpha}{b-1}+\left[0, \frac{t(N-1)}{b^k(b-1)}\right]$ contains at most $\left\lceil \frac{t(N-1)}{b^k(b-1)}\right\rceil$ integers for each $I\in\{0,1,\ldots,N-1\}^k$.
For $k\ge 1$,
\begin{eqnarray*}
	\#(T(b, t\CC)\cap\Z\setminus\{0\})&\le& N^k\left\lceil \frac{t(N-1)}{b^k(b-1)}\right\rceil-1\\
	&=&N^k\left(\frac{t(N-1)}{b^k(b-1)}-\left\{\frac{t(N-1)}{b^k(b-1)}\right\}+1\right)-1\\
	&<&\frac{t(N-1)}{q^k(b-1)}+N^k.
\end{eqnarray*}
\eproof

\begin{lem}\label{lem3.2}
	Let $t\in\N^\ast$ with $\gcd(t,N)=1$ be a primitive incomplete number of the spectral pair $(\mu_{b, qD}, \Lambda(b, \CC))$. If $\{x_k\}_{k=0}^\infty$ is a non-zero integer periodic sequence w.r.t. $(b, t\CC)$, then $O_b(t)$ is its minimal positive period and
	\[\left\{x_1,\ldots,x_{O_b(t)}\right\}=\left\{b^ix_k\pmod t:1\le i\le O_b(t)\right\}, \quad \text{for any $k\ge 1$.}\]
\end{lem}
\begin{proof}
	By Corollary \ref{coro2.9} we have $\gcd(t, b)=1$. Suppose that $n$ is the minimal positive period of $\{x_k\}_{k=0}^{\infty}$. Then
	\begin{equation}\label{xklies}
		x_k\in T(b, t\CC)\subset \left( -\frac{tN}{2(b-1)}, \frac{tN}{2(b-1)}\right]\subset \left(-\frac t2, \frac t2\right)
	\end{equation}
	for all $k$. From $bx_1=x_0+tc_0$ for some $c_0\in\CC$ we have $x_1\equiv b^{O_b(t)-1}x_0\pmod t$. Then $x_{O_b(t)}\equiv x_0\pmod t$ and thus $x_{O_b(t)}=x_0$ by \eqref{xklies}. Corollary \ref{coro2.12} indicates that $n\mid O_b(t)$. If $n<O_b(t)$, then similarly we have $x_0=x_n\equiv b^{O_b(t)-n}x_0\pmod t$, that is, there exists $k_0\in\N$ such that
	\[b^{O_b(t)-n}x_0=x_0+tk_0.\]
	By Corollary \ref{coro2.9} again we obtain that $b^{O_b(t)-n}\equiv1\pmod t$, which contradicts the definition of $O_b(t)$.
\end{proof}

\begin{coro}\label{coro4.10}
	Let $p$ be a prime such that
	\[O_b(p)\ge \frac{p(N-1)}{q^k(b-1)}+N^k\qquad \text{for some $k\ge 1$.}\]
	Then $p$ is a complete number of  of the spectral pair $(\mu_{b, q\D}, \La(b, \CC))$.
\end{coro}
\proof Suppose on the contrary that $p$ is an incomplete number. Then $p$ is primitive since $p$ is a prime number. Corollary \ref{coro2.11} and Lemma \ref{lem3.2} imply that $T(b,t\CC)\cap\Z\setminus\{0\}$ has cardinality no less than $O_b(t)$, which contradicts Lemma \ref{lem4.9}.
\eproof

Now we use the following result, which was proved by Erd\"os and Murty \cite[Theorem 1]{EM1999}.
\begin{theo}\label{theoerdo}
	Let $a>1$ and let $\epsilon(p)$ be any function tending to zero as $p\to\infty$. For all but $o(x/\log x)$ primes $p\le x$,
	\begin{eqnarray*}
		O_a(p)\ge p^{\frac 1{2+\epsilon(p)}}.
	\end{eqnarray*}
\end{theo}
According to the prime number theorem, the above result says that
\begin{eqnarray}\label{eq4.3}
	&&\lim_{ x\to\infty}\frac{\#\{p\le x\text{ is a prime}: O_a(p)\ge p^{\frac 1{2+\epsilon(p)}}\}}{\#\{p\text{ is a prime}: p\le x\}}\nonumber\\
	&=&\lim_{x\to\infty}\frac{\#\{p\le x\text{ is a prime}: O_a(p)\ge p^{\frac 1{2+\epsilon(p)}}\}}{x/\log x}=1.
\end{eqnarray}
Let $\Gamma\subset\N$. For $\alpha>0$, the {\it upper and lower Banach density } is defined by
\begin{eqnarray*}
	D_\alpha^+(\Gamma)=\limsup_{h\to\infty}\frac{\#(\Gamma\cap [0, h])}{h^\alpha} \quad\text{and}\quad D_\alpha^-(\Gamma)=\liminf_{h\to\infty}\frac{\#(\Gamma\cap [0, h])}{h^\alpha}
\end{eqnarray*}
respectively.
\begin{prop}\label{firddst}
	Let $\Gamma$ be the set of all positive complete numbers. If $q>N$, then
	\begin{eqnarray*}
		\liminf_{x\to\infty}\frac{\#\Gamma\cap[0,x]}{x/\log x}\ge 1.
	\end{eqnarray*}
\end{prop}
\proof
We take $\epsilon(p)\le\log_pq$. Then
\begin{eqnarray*}
	&&\log_N\frac12 p^{\frac 1{2+\epsilon(p)}}-\log_q 2 p^{\frac {1+\epsilon(p)}{2+\epsilon(p)}}\\
	&=&\log_N\frac12+{\frac 1{2+\epsilon(p)}}\log_N p-\log_q 2-{\frac {1+\epsilon(p)}{2+\epsilon(p)}}\log_q p\\
	&\ge&-\log_N2-\log_q2+{\frac 1{2+\epsilon(p)}}(\log_N p-\log_qp)-1,
\end{eqnarray*}
which is larger than 1 if $q>N$ and $p\ge p_0$ for some $p_0$ depending only on both $N$ and $q$. Then there exists $k_0$ such that
\[\log_q 2 p^{\frac {1+\epsilon(p)}{2+\epsilon(p)}}\le k_0\le \log_N\frac12 p^{\frac 1{2+\epsilon(p)}}.\,\]
which is equivalent to that
\begin{eqnarray*}
	\frac p{q^{k_0}}\le \frac12p^{\frac 1{2+\epsilon(p)}}\quad \text{and}\quad N^{k_0}\le \frac12 p^{\frac 1{2+\epsilon(p)}}.
\end{eqnarray*}
Based on the above analysis
\[p^{\frac 1{2+\epsilon(p)}}\ge\frac p{q^{k-0}}+N^{k_0}>\frac{p(N-1)}{q^{k_0}(b-1)}+N^{k_0}.\]

For $x>0$ we denote
\[\Gamma^*(x)=\{p\text{ is a prime}:p_0\le p\le x \text{ and } O_b(p)\ge p^{\frac 1{2+\epsilon(p)}}\}.\]
By \eqref{eq4.3} we have
\begin{eqnarray*}
	\lim_{x\to\infty}\frac{\Gamma^*(x)}{x/\log x}=1.
\end{eqnarray*}
Corollary \ref{coro4.10} implies that $\Gamma^\ast\subset\Gamma\cap[0,x]$.
This complete the proof.
\eproof
\begin{lem}\label{lem4.13}
	Denote $\alpha_k=1+b+\cdots+b^{k-1}$. Then $\alpha_k$ is an incomplete number for each $k\ge b$. Moreover, $\gcd(\alpha_k, \alpha_l)=1$ if and only if $\gcd(k, l)=1$.
\end{lem}
\proof According to (iii) of Theorem \ref{theo2.5} and
\begin{equation*}
	\alpha_k\frac{b+\cdots+b^{b-1}}{b^{k}-1}=\frac{b-1+b-1+b^2-1+\cdots+b^{b-1}-1}{b-1}\in \alpha_kT(b, \CC)\cap\N,
\end{equation*}
$\alpha_k$ is an incomplete number for each $k\ge b$.

For $k<l$ and write $l=sk+t$, $0<t<k$, then
$$\gcd(\alpha_k, \alpha_l)=\gcd(\alpha_k, b^{sk}+b^{sk+1}+\cdots+b^{sk+t-1})=\gcd(\alpha_k, \alpha_t).$$
This implies the desired result. \eproof

The following proposition gives an upper bound of $D_1^+(\Gamma)$.
\begin{prop}\label{secddst}
	Let $\Gamma\subset\N^\ast$ be the set of all positive complete numbers. Then $$\limsup_{x\to\infty}\frac{\#(\Gamma\cap[0,x])}{x}\le c\cdot\frac{\varphi(N)}{N},$$
	where $\varphi$ is Euler function and $0<c<1$ is a constant depending on $b$.
\end{prop}
\begin{proof}
	Denote by $N_1,N_2,\ldots,N_s$ all distinct prime factors of $N$, and let $\{p_i\}_{i=1}^{\infty}$ be the increasing sequence of prime numbers greater than $b$. For convenience, we define $N_{s+i}:=\alpha_{p_i}$, where $\alpha_{p_i}$ is given as in Lemma \ref{lem4.13}. By Lemma \ref{lem4.13}, we conclude that $N_1,\ldots,N_{s+n}$ are mutually coprime and are incomplete numbers.
	From Corollary \ref{coronn}, for any positive integer $n$, it is obvious that
	$$\Gamma\subset\N^\ast\setminus\bigg(\bigcup_{i=1}^{s+n}N_i\N^\ast\bigg).$$
	We let $A_i=N_i\N^\ast\cap[0,x]$. The Inclusion-Exclusion Principle indicates that
	\begin{eqnarray*}
			\#\Bigg( \bigcup_{i=1}^{s+n} A_i\Bigg) &=& \sum_{i=1}^{n} \#(A_i) - \sum_{1 \le i < j \le s+n} \#(A_i \cap A_j) + \sum_{1 \le i < j < k \le n} \#(A_i \cap A_j \cap A_k) \\
			&&\qquad-\cdots + (-1)^{s+n+1} \#(A_1 \cap A_2 \cap \cdots \cap A_n)\\
			&=& \sum_{i=1}^{n} \left\lfloor\frac{x}{N_i}\right\rfloor - \sum_{1 \le i < j \le s+n} \left\lfloor\frac{x}{N_i N_j}\right\rfloor + \sum_{1 \le i < j < k \le n} \left\lfloor\frac{x}{N_i N_j N_k}\right\rfloor \\
			&&\qquad-\cdots + (-1)^{s+n+1} \left\lfloor\frac{x}{N_1 N_2 \cdots N_{s+n}}\right\rfloor\\
			&=&\sum_{i=1}^{n}\frac{x}{N_i} - \sum_{1 \le i < j \le s+n} \frac{x}{N_i N_j} + \sum_{1 \le i < j < k \le n} \frac{x}{N_i N_j N_k} \\
			&&\qquad-\cdots + (-1)^{s+n+1}\frac{x}{N_1 N_2 \cdots N_{s+n}}+C_n
	\end{eqnarray*}
	for some finite number $C_n$. Then
	\begin{eqnarray*}
		\limsup_{x\to\infty}\frac{\#(\Gamma\cap[0,x])}{x}&\le&1-\sum_{i=1}^{n}\frac{1}{N_i} - \sum_{1 \le i < j \le s+n} \frac{1}{N_i N_j} + \sum_{1 \le i < j < k \le n} \frac{1}{N_i N_j N_k} \\
		&&\qquad-\cdots + (-1)^{s+n}\frac{1}{N_1 N_2 \cdots N_{s+n}}\\
		&=&\prod_{i=1}^{s+n}\left( 1-\frac{1}{N_i}\right) \\
		&\le&\frac{\varphi(N)}{N}\prod_{i=s+1}^n\left( 1-\frac{1}{N_i}\right).
	\end{eqnarray*}
	Letting $n\rightarrow\infty$, we have
	$$\limsup_{x\to\infty}\frac{\#(\Gamma\cap[0,x])}{x}\le \frac{\varphi(N)}{N}\cdot \prod_{p\ge b, prime}\left(1-\frac{1}{\alpha_{p}} \right) .$$
\end{proof}

{\noindent\bf Proof Theorem \ref{thmdeig}.} They can be directly obtained from Proposition \ref{firddst} and Proposition \ref{secddst}.

\section{Conjecture \ref{conj1} for the spectral pair $(\mu_{b, q\D}, \Lambda(b, \CC))$}\label{secprconj}

In this section, we mainly focus on the positive integers that are coprime with $b$. The most representative ones are primitive numbers (complete numbers and incomplete numbers).

\begin{theo}\label{theo1.10}
	Let $t\in\N^\ast$ with $\gcd(t,N)=1$ be a primitive incomplete number. Then the following statements hold:
	\begin{enumerate}
		\item [(i)] $\gcd(t, b)=1$.
		\item [(ii)] All non-zero integer periodic sequence w.r.t. $(b, t\CC)$ have minimal positive period $O_b(t)$.
		\item [(iii)] There exist $c_0, c_1, \ldots, c_{O_b(t)-1}\in \CC$, not all equal to $0$, such that
		\begin{eqnarray*}
			t\frac{c_0+bc_1+\cdots+b^{O_b(t)-1}c_{O_b(t)-1}}{b^{n}-1}\in T(b, t\CC)\cap\Z\setminus\{0\}.
		\end{eqnarray*}
		\item [(iv)] Let $\{x_k\}_{k=0}^{\infty}$ be a non-zero integer periodic sequence w.r.t. $(b, t\CC)$. Then
		\[\left\{x_1,\ldots,x_{O_b(t)}\right\}=\left\{b^ix_k\pmod t:1\le i\le O_b(t)\right\}\subset T(b, t\CC)\cap \Z\setminus\{0\}, \]
		for any $k\ge 1$.
	\end{enumerate}
	Moreover, the statements $(ii)$, $(iii)$ and $(iv)$ are equivalent and guarantee the property of incomplete number of $t$.
\end{theo}

\begin{proof}
	Theorem \ref{theo2.5}, Corollary \ref{coro2.9} and Lemma \ref{lem3.2} indicate the proof of this theorem.
\end{proof}
{\noindent\bf Remark.} Let $t \in \N^\ast$ with $\gcd(t, b) = 1$. There exist cases where $t$ is not a primitive incomplete number, but all non-zero integer periodic sequences w.r.t. $(b, t\CC)$ have a minimal positive period of $O_b(t)$. An illustrative example of this can be found in Example \ref{exam3.7}.

\begin{exam}\label{exam3.7}
	For the case $N=3$, $q=2$. Let $t=72559411$. Then $t$ is an incomplete number and all non-zero integer periodic sequence w.r.t. $(6,t\{-1,0,1\})$ have minimal positive period $O_6(t)$, but $t$ is not primitive.
	\begin{proof}
		Note that $72559411=23\times 3154757$ and $O_6(23)=O_6(3154757)=O_6(72559411)=11$. Example \ref{exam2.12} implies that the prime $3154757$ is an incomplete number, hence primitive. By Corollary \ref{coronn} and the definition of primitive incomplete number, $t$ is an incomplete number and not primitive. Suppose that $\{x_k\}_{k=0}^\infty$ is a non-zero integer periodic sequence w.r.t. $(6,t\{-1,0,1\})$ and its minimal positive period is $n$. Then $x_0=6x_1\pmod t$. By induction,
		$$x_{O_6(t)}\equiv 6^{O_6(t)}x_{O_6(t)}\equiv x_0\pmod t.$$
		Since all $x_k$ lie in $T(6,t\{-1,0,1\})\subset\left( -\frac{t}{2},\frac{t}{2}\right)$, we have $x_{O_6(t)}=x_0$. This forces $n\mid O_6(t)$. It is easy to check that $n$ cannot be $1$. Hence $n=O_6(t)=11$.
	\end{proof}
\end{exam}
The following proposition is similar to the one given by Dutkay and Haussermann \cite[Theorem 2.3]{DH16}. We give its proof here for completeness.
\begin{prop}\label{prop4.2}
	There are infinitely many primitive incomplete numbers.
\end{prop}
\begin{proof}
	Noting that $b-1$ is an incomplete number, there exists at least one primitive incomplete number that divides $b-1$. We choose an arbitrary finite set of primitive incomplete numbers, where those that are not factors of $b-1$ are denoted by $t_1,t_2,\ldots,t_k$. Now we construct a new primitive incomplete number that is not in the finite set we choose.
	
	We denote $n:=(b-1)O_b((b-1)^2)O_b(t_1)\cdots O_b(t_k)\ge(b-1)$. Then $b^{n+1}-1\equiv b-1\pmod{(b-1)^2,t_1,\ldots,t_k}$. Letting
	$$t:=\frac{b^{n+1}-1}{b-1}\in\N,$$
	we claim that $t$ is not divisible by any one of $b-1$, $t_1$, $\ldots$, $t_k$. Otherwise, $b^{n+1}-1$ is divisible by one of $(b-1)^2$, $t_1$, $\ldots$, $t_k$, this contradicts the fact $b^{n+1}-1\equiv b-1\pmod{(b-1)^2,t_1,\ldots,t_k}$.
	
	Clearly,
	$$t\frac{1+b+\cdots+b^{b-2}}{b^{n+1}-1}=\frac{1+b+\cdots+b^{b-2}}{b-1}\in T(b,t\CC)\cap\Z\setminus\{0\}.$$
	By Theorem \ref{theo2.5}, $t$ is an incomplete number. We let $t'\mid t$ be a primitive incomplete number. We conclude that $t'\neq t_i$ for any $1\le i\le k$ by our claim and then  $t'\mid (b-1)$. Since $O_b((b-1)^2)\mid n$, there exists $l\in\Z$ such that
	$$t(b-1)=b^{n+1}-1=b-1+l(b-1)^2,$$
	and then $t=1+l(b-1)$. This forces $t'=1$, which contradicts the fact $t'$ is a primitive incomplete number. Hence, we complete the proof.
\end{proof}

\begin{theo}\label{theo4.1}
	Let $t\in\N^\ast$ with $\gcd(t,b)=1$ and cyclic group $G(b,t)=\{b^l\pmod t:1\le l\le O_b(t)\}$.
	\begin{enumerate}
		\item [(i)] If $q\ge 3$, $2\le |a|\le q-1$ and $a\in G(b,t)$, then $t$ is a complete number.
		\item[(ii)] If $q=2$, $|a|=2$ and $a\in G(b,t)$, then $t$ is a complete number.
	\end{enumerate}
\end{theo}
\proof
Suppose by contradiction that $t$ is an incomplete number, without loss of generality, we can assume that $t$ is a primitive incomplete number. Then there exists a non-zero integer periodic sequence $\{x_k\}_{k=0}^\infty$ w.r.t. $(b,t\CC)$. For $k\in\N$,
\begin{eqnarray}\label{eq4.2}
	x_k\in T(b,t\CC)\cap\Z\subset\left\{
	\begin{array}{ll}
		\bigg[ -\frac{t(N-2)}{2(b-1)}, \frac{tN}{2(b-1)}\bigg]\cap\Z , & \hbox{if $N$ is even;} \\
		\bigg[ -\frac{t(N-1)}{2(b-1)}, \frac{t(N-1)}{2(b-1)}\bigg]\cap\Z, & \hbox{if $N$ is odd.}
	\end{array}
	\right.
\end{eqnarray}
We let $x_0$ be the point which is maximum in modulus. Then $ax_0\pmod t=x_{k_0}$ for some $k_0\in\N$ by the assumption that there exists $l<O_b(t)$ such that $b^{l}\equiv a\pmod t$ and Lemma \ref{lem3.2}. According to \eqref{eq4.2} and $|a|\le q-1<\frac{b-1}m$ one has that $|ax_0|<\frac t2$ and thus $ax_0=x_{k_0}$. Then $|x_{k_0}|=|ax_0|>|x_0|$, which contradicts $|x_0|=\max_{k\in\N}|x_k|$.

When $b=2N$ for some $N>2$ and $|a|=2$. If $N$ is odd, then similarly, by \eqref{eq4.2} we have
\[|ax_0|\le \frac{t(N-1)}{b-1}<\frac t2,\]
which yields that $ax_0=x_{k_0}$ and contradicts the maximality of $|x_0|$. If $N>2$ is even and $ax_0\in(-\frac t2, \frac t2]$, then it is impossible by the same idea.

The remaining case is that $N>2$ is even and $2|x_0|>\frac t2$ or $ax_0=-\frac t2$. Therefore by \eqref{eq4.2} again, $|ax_0|=2|x_0|<t$ and
\begin{eqnarray*}
	x_{k_0}=ax_0\pmod t=\left\{
	\begin{array}{ll}
		ax_0+t, & \text{if }ax_0<0; \\
		ax_0-t, & \text{if }ax_0>0.
	\end{array}
	\right.
\end{eqnarray*}
Then
$$t\le\vert ax_0\vert+\vert x_{k_0}\vert\le\frac{3tN}{2(b-1)}<t,$$
which yields a contradiction.\ez

\begin{coro}\label{coro4.4}
	If $d\mid q$ for some $1<d<q$, then $dNb^k\pm1$ is a complete number for any $k\in\N^\ast$.
\end{coro}
\begin{proof}
	Writing $t=dNb^k-1$ for any chosen $k\in\N^\ast$, we obtain $b^{k+1}-\frac{q}{d}t=\frac{q}{d}\in\Z$ and $2\le\frac{q}{d}\le q-1$. So $dmb^k-1$ is a complete number and similarly for $dmb^k+1$.
\end{proof}

The cases $\mu_{4,\{0,2\}}$ and $\mu_{2q, \{0, q\}}$ of Theorem \ref{theoinfinite} have been proved in \cite[Theorem 2.3, Corollary 2.9]{DH16} and \cite[Corollary 2.12, Theorem 2.17]{DK18}, respectively. Now we prove the Theorem in our setting:

\noindent{\bf Proof of Theorem \ref{theoinfinite}.} According to Proposition \ref{prop4.2}, we only need to prove that there are infinite primitive complete numbers.

Case 1. $q$ is not a prime. For $k\ge 1$, let $t_k=dNb^k+1$ for some $d\mid q$ and $1<d<q$, then $\gcd(t_k, b)=1$ and all $t_k$ are complete numbers for $k\ge 1$ by Corollary \ref{coro4.4}.

Case 2. $q$ is  an odd prime and $m$ is odd. Note that $\gcd(b^k+2,b)=1$ for $k\geq1$, by the following simple fact: $\frac{b+1}{2}(b^k+2)+(-\frac{b^{k}+2+b^{k-1}}{2})\cdot b=1$ since $b=Nq$ is odd. Then all $t_k:=b^k+ 2$ are primitive complete numbers by Theorem \ref{theo4.1}.

Case 3. $q$ is  an odd prime and $N$ is even. Denote $t_k=\frac N2qb^k+1$ for $k\ge 1$. Then $\gcd(t_k, b)=1$ and $b^{k+1}\equiv -2{\pmod {t_k}}$. Then all $t_k$ are primitive complete numbers by Theorem \ref{theo4.1}.

Case 4. $q=2$. Denote $t_k=Nb^k+1$ for $k\ge 1$. Then $\gcd(t_k, b)=1$ and $b^{k+1}\equiv -2\pmod {t_k}$ and the desired result follows from Theorem \ref{theo4.1} again. \ez

\section{Sufficient conditions for the complete numbers}\label{secsuffi}

Dutkay et al. \cite{DH16} first discovered the relationship between (incomplete) complete numbers and number theory. Here we introduce their related results and some new results.

\begin{lem}\label{lem5.1}
	Let $t$ be an integer with $\gcd(t,b)=1$. Then $t$ is a divisor of $b^{n}-1$ if and only if $O_{b}(t)\mid n$.
\end{lem}
\begin{proof}
	Suppose that $O_{b}(t)\mid n$. Then  $b^n\equiv 1\pmod {t}$, i.e., $t$ is a divisor of $b^{n}-1$.
	Conversely, suppose that $t$ is a divisor of $b^{n}-1$. If $O_b(t)\nmid n$, then we denote $n=kO_b(t)+m$, where $k\in\N$ and $1\le m\le O_b(t)-1$. In this case
	$$1\equiv b^n=b^{kO_b(t)+m}\equiv b^m\pmod t.$$
	This contradicts the definition of $O_b(t)$.
\end{proof}
Following \cite{DH16}, we define
\begin{defi}\label{defi4.3}
	For a prime number $p$, we denote by $\ell_b(p)$ the largest positive integer $k$ such that $O_b(p^k)=O_b(p)$. The prime $p$ is referred to as a simple prime if $\ell_b(p)=1$.
\end{defi}
The following three results were proved by Dutkay and Kraus in Section 2.4 of \cite{DK18}.
\begin{prop}\label{prop4.4}
	Let $s$ and $t$ be relatively prime numbers that $\gcd(st,b)=1$. Then
	$$O_b(st)=\operatorname{lcm}\big(O_b(s),O_b(t)\big).$$
\end{prop}

\begin{prop}\label{prop5.4}
	Let $p$ be prime that is coprime with $b$. Then $O_b(p^k)=O_b(p)$ for $k\le\ell_b(p)$ and $O_b(p^k)=p^{k-\ell_b(p)}O_b(p)$ for all $k>\ell_b(p)$.
\end{prop}

\begin{prop}\label{prop4.6}
	Suppose $p_1,\ldots,p_N$ are distinct primes with $\gcd(p_1\cdots p_N,b)=1$ and $k_1,\ldots,k_N\ge 0$. For $1\le i\le N$, let $l_i$ be the largest integer such that $p^{l_i}_i$ divides $\operatorname{lcm}(O_b(p_1),\ldots,O_b(p_N))$ . Then
	\begin{equation*}
		O_b({p_1^{k_1}\cdots p_N^{k_N}})=\Big(\prod_{i=1}^{N}p_i^{\max\{k_i-\ell_b(p_i)-l_i,0\}}\Big)\operatorname{lcm}\big( O_b(p_1),\ldots,O_b(p_N)\big).
	\end{equation*}
\end{prop}
As usual $\lfloor x \rfloor=\max\{y\in\Z:y\le x\}$ and $\lceil x \rceil=\min\{y\in\Z:y\geq x\}$ for $x\in\R$.
\begin{lem}\label{lem5.6}
	Let $s,t>1$ with $\gcd(st,b)=1$ be two complete numbers of the spectral pair $(\mu_{b, q\D}, \La(b, \CC))$. Assume that
	$$\frac{O_{b}(st)}{O_{b}(t)}\ge\frac{2b^2-2Nb+b+2N-3+Ns}{q(b-1)},$$
	then $st$ is not a primitive incomplete number.
\end{lem}
\begin{proof}
	Suppose on the contrary, $st$ is a primitive incomplete number. Then there exists a non-zero integer periodic sequence $\{x_k\}_{k=0}^\infty$ w.r.t. $(b,st\CC)$.
	Let $S=\{x_1,\ldots,x_{O_b(t)}\}$, by Lemma \ref{lem3.2} one has
	\begin{eqnarray*}
		S&=&\{b^ix_0\pmod {st}: 1\le i\le O_{b}(st)\}\\
		&\subset& T(b, st\CC)\cap\Z\setminus\{0\}\subset \left(-\frac{stN}{2(b-1)}, \frac{stN}{2(b-1)}\right].
	\end{eqnarray*}
	We define an equivalent class in $S$ by $x_i\sim x_j$ if $[x_i]_t=[x_j]_t$ for $x_i,x_j\in S$. Here and in the following sections, $[a]_t$ means the congruence class of $a$ modulo $t$. Then $b^ix_0\pmod{st}\sim b^jx_0\pmod{st}$ if and only if $t\mid b^ix_0(b^{j-i}-1)$ for $1\le i<j\le O_{b}(st)$. Note that $\gcd(t, x_0)=1$, as stated in Corollary \ref{coro2.9}. Then, by applying Lemma \ref{lem5.1},
	$$b^ix_0\pmod{st}\sim b^jx_0\pmod{st}\Longleftrightarrow O_{b}(t)\mid(j-i).$$
	Let $M:=\frac{O_{b}(st)}{O_{b}(t)}\in\N$. Then all equivalent classes in $S$
	$$\left[b^ix_0\pmod{st}\right]_t=\left\{b^{i+kO_{b}(t)}x_0\pmod{st}:0\le k\le M-1\right\}$$
	have the same cardinality $M$, i.e., $S/{\sim}=\left\{\left[b^ix_0\pmod{st}\right]_t\right\}_{i=1}^{O_b(t)}$.
	
	We take an $x_i$ in $S$ and $y_i\equiv x_i\pmod t\in (-\frac t2, \frac t2]$. Then there exist $M$ values of $k$ such that $y_i+kt$ in the congruence class $[x_i]_t$ in $S$.	Noticing that any element in $S$ is congruent to $-stc$ modulo $b$ for some $c\in\CC$. Then we have for $M$ values of $k$, $y_i+kt\equiv -stc_k\pmod b$ for some $c_k\in \CC$. Then $k\equiv t^{O_t(b)-1}(-y_i-stc_k)\pmod b$ and we consider the possible values of $k$. Since there are $M$ values of $k$, at least one of the two intervals $(-\infty,0]$ and $[0,\infty )$ lies $\lceil\frac{M}{2}\rceil$ values of $k$. Without loss of generality, we suppose that $\lceil\frac{M}{2}\rceil$ values of $k$ in $[0,\infty )$.
	
	We consider the case that there are $\lceil\frac{M}{2}\rceil$ values of $k$ in the interval $[0,\infty )$. $k\equiv t^{O_t(b)-1}(-y_i-stc_k)\pmod b$ implies that for each $n\ge 0$, the set $K_n:=\{nb,nb+1,\ldots,nb+b-1\}$ contains at most $N$ values of $k$. If $\lceil\frac{M}{2}\rceil$ is divisible by $N$, then there are at most ${\lceil\frac{M}{2}\rceil}-N$ values of $k$ belonging to $K_0\cup K_1\cup\cdots\cup K_{m_1-1}$, where $m_1={\frac{\lceil\frac{M}{2}\rceil}{N}-1}$. Therefore, the largest value of $k$ satisfies that
	$$k\ge m_1b+N-1=q{\bigg\lceil\frac{M}{2}\bigg\rceil}-b+N-1.$$
	If $\lceil\frac{M}{2}\rceil$ is not divisible by $m$, then we suppose that $\lceil\frac{M}{2}\rceil=m_2N+\alpha$, $1\le \alpha\le N-1$. There are at most $m_2N$ values of $k$ belong to the set $K_0\cup K_1\cup\cdots\cup K_{m_2-1}$. Then the largest value of $k$ satisfies that
	$$k\ge m_2b-1+\alpha=\frac{\big\lceil\frac{M}{2}\big\rceil-\alpha}{N}b-1+\alpha=q\bigg\lceil\frac{M}{2}\bigg\rceil-(q-1)\alpha-1\geq q\bigg\lceil\frac{M}{2}\bigg\rceil-b+N+q-2$$
	by $\alpha\le N-1$. Anyway, $k\ge \frac{qM}{2}-b+N-1$.
	Thus, there is an element in $S$,
	\begin{equation*}
		y_i+kt>\left( q{\frac{M}{2}}-b+N-\frac32\right) t\ge\frac{stN}{2(b-1)},
	\end{equation*}
	which contradicts the range of $S$. \end{proof}

Dutkay and Kraus proved the case $N=2$ of the following Theorem \cite[Theorem 2.42]{DK18}. We focus on the case $N\ge 3$.
\begin{theo}\label{theodut}
	Let $p_1,\ldots,p_m\ge 2N+1$ be distinct primes with $\gcd(p_1\cdots p_m,b)=1$. For $1\le i\le m$, let $l_i$ be the largest integer such that $p_i^{l_i}\mid\operatorname{lcm}\big( O_b(p_1),\ldots,O_b(p_m)\big)$. Assume that the integer $p_1^{\ell_b(p_1)+l_1}\cdots p_m^{\ell_b(p_m)+l_m}$ is a complete number. Then so is $p_1^{k_1}\cdots p_m^{k_m}$  for any $k_1,\ldots,k_m\ge 0$.
\end{theo}
\proof
Suppose that there exist some nonnegative integers $k_1,\ldots,k_m$ such that $p_1^{k_1}\cdots p_m^{k_m}$ is an incomplete number. Then there is a primitive incomplete number that divides $p_1^{k_1}\cdots p_m^{k_m}$. Without loss of generality, we denote such proper divisor by $p_1^{k_1}\cdots p_m^{k_m}$. The hypothesis implies that there exists $k_j\ge\ell_b(p_j)+l_j+1$ for some $1\le j\le m$. Relabeling the primes $p_i$, we assume $k_1\ge\ell_b(p_1)+l_1+1$. Proposition \ref{prop4.6} implies that
$$O_b(p_1^{k_1}\cdots p_m^{k_m})=p_1^{k_1-\ell_b(p_1)-l_1}O_b(p_1^{\ell_b(p_1)+l_1}p_2^{k_2}\cdots p_m^{k_m}).$$
Write $s=p_1^{k_1-\ell_b(p_1)-l_1}$ and $t=p_1^{\ell_b(p_1)+l_1}p_2^{k_2}\cdots p_m^{k_m}$.	We claim that if $p_1>2N$
\begin{equation}\nonumber
	\begin{split}
		\frac{O_b(st)}{O_b(t)}=s\ge\frac{2b^2-2Nb+b+2N-3+Ns}{q(b-1)}.
	\end{split}
\end{equation}
In fact,
$$
s=\frac{Ns}{q(b-1)}+\left( 1-\frac{N}{q(b-1)}\right)s\ge\frac{Ns}{q(b-1)}+\left( 1-\frac{N}{q(b-1)}\right)(2N+1),\\
$$
we need to check that
$$\left( 1-\frac{N}{q(b-1)}\right)(2N+1)\ge\frac{2b^2-2Nb+b+2N-3}{q(b-1)},$$
which is equivalent to that
$$2N+q-3\ge\frac{2N+1}{q-\frac{1}{N}}$$
by calculation. The above inequality holds since we suppose $m\ge3$ and $q\ge2$. Indeed,
$$\frac{2N+1}{q-\frac{1}{N}}\le\frac{2N+1}{2-\frac13}<2N-1\le2N+q-3.$$

By Lemma \ref{lem5.6}, $p_1^{k_1}\cdots p_m^{k_m}$ is not a primitive incomplete number, a contradiction.\ez

{\noindent\bf Remark.} In this theorem and the following corollary, when $N=2$, the condition that relevant primes are greater than $2N+1$ is not required.
\begin{coro}\label{coro5.8}
	Let $p$ be a prime with $p\ge 2N+1$ and $\gcd(p, b)=1$. If $p^{\ell_b(p)}$ is a complete number of the spectral pair $(\mu_{b, q\D}, \La(b, \CC))$, then all $p^k$ are complete numbers.
\end{coro}

\section{The simplest case ($N=2$)}\label{sectwo}

Now we give the following counterexamples to illustrate that the answer to question 1.3 is negative for some cases.
\begin{exam}\label{exam2.12}
	For the case $N=3$, $q=2$ and let $p=3154757$. Using the computer numerical calculations, $p$ is a prime and there exists a non-zero integer periodic sequence w.r.t. $(6,p\{-1,0,1\})$. We list the terms of a period in order:
	\begin{equation}\nonumber
		\begin{split}
			\{2015,-525457,-613369,&-628021,-630463,-630870,\\
			&-105145,-543317,435240,72540,12090\}.
		\end{split}
	\end{equation}
	Then the prime $p$ is an incomplete number of the spectral pair $(\mu_{6,\{0,1,2\}},\Lambda(6,\{-1,0,1\}))$.
\end{exam}
The following example was given by \cite[Remark 2.14]{DK18}.
\begin{exam}\label{exam2.13}
	For the case $N=2$, $q=3$ and let $p=55987$. There exists a non-zero integer periodic sequence w.r.t. $(6,p\{0,1\})$. Terms of a period in order are
	$$\{311,9383,10895,11147,11189,11196,1866\}$$
	and then the prime $p$ is an incomplete number of the spectral pair $(\mu_{6,\{0,3\}},\Lambda(6,\{0,1\}))$.
\end{exam}

In this section we study the case $N=2$. The special case $b=4$ has been studied by  Dutkay and  Haussermann \cite{DH16}. Here we study the following question in situation $b=2^r$ for  $r>2$:

{\it Let $p$ be an odd prime. Are $p^n$ for $n\ge 1$ complete numbers of the spectral pair $(\mu_{2^r, \{0, 2^{r-1}\}}, \Lambda(2^r, \{0, 1\})$?}

{\noindent The answer to this question is positive for $r=2$ and $p>3$, which was proved by \cite[Theorem 2.10]{DH16} with different methods.}

In Section 2, we have proved that $2^r-1$ is incomplete number and the premise for $t$ being a complete number of the spectral pair $(\mu_{2^r, \{0, 2^{r-1}\}}, \La(2^r, \{0, 1\})$ is $\gcd(t,2)=1$. So odd prime numbers that are not equal to $2^r-1$ are concerned in this section. We start with the following known lemma,  For completeness, we give its proof.
\begin{lem}\label{lem6.1}
	Let $\gcd(t, b)=1$. Then for $r\ge 1$
	\begin{eqnarray*}
		O_{b^r}(t)=\frac{O_b(t)}{\gcd(O_b(t), r)}.
	\end{eqnarray*}
\end{lem}
\proof Denote $d=\gcd(O_b(t), r)$ and write $O_b(t)=d\al$, $r=dr'$. By the definition of the order we have $b^{rO_{b^r}(t)}\equiv 1\pmod t$. Then $O_b(t)\mid rO_{b^r}(t)$. This implies that $\al\mid O_{b^r}(t)$, i.e,
\[\frac{O_b(t)}{\gcd(O_b(t), r)}\mid  O_{b^r}(t).\]
Conversely, from $b^{O_b(t)}\equiv1\pmod t$ one has
$$b^{r'O_b(t)}=b^{r\al}\equiv1\pmod t,$$
which yields that $O_{b^r}(t)\mid\al$. Then the converse holds. \ez

Using Lemma \ref{lem6.1}, we can prove the  following several lemmas:
\begin{lem}\label{lem6.2}
	Let $t>1$ be an odd number. If there exists $i\ge 0$ such that
	$$2^{ir}\equiv a\pmod t,\qquad \text{for some $k=-1$ or $2\le|a|\le 2^{r-1}$}$$
	then $t$ is a complete number of the spectral pair $(\mu_{2^r, \{0, 2^{r-1}\}}, \Lambda(2^r, \{0, 1\}))$.
\end{lem}
\proof  When  $|a|\le \max\{2^{r-1}-1, 2\}$, the assertion has been proved by Theorem \ref{theo4.1} in this paper. When $a=-1$, with loss of generality  we can assume that $t$ is a primitive incomplete number. Then there exists a non-zero integer periodic sequence $\{x_k\}_{k=0}^\infty$ w.r.t. $(2^r, \{0,  t\})$ and for any $k\in\N$,
\begin{eqnarray*}
	x_k\in T(2^r,\{0, t\})\setminus\{0\}\cap\Z\subset \left(0, \frac t{2^r-1}\right].
\end{eqnarray*}
Let $x_0=\max_{k\in\N} x_k$. Then $1\le x_0\le \frac t{2^r-1}$. By the assumption and Lemma \ref{lem3.2} we have
$$2^{ir}x_0\pmod t=-x_0\pmod t=t-x_0.$$
However, $t-x_0\ge t-\frac t{2^r-1}>\frac t{2^r-1}$, which yields a contradiction.

When $a=2^{r-1}$, similarly
$$2^{ir}x_0\pmod t=2^{r-1}x_0\pmod t.$$
If $2^{r-1}x_0\le \frac t2$, then $2^{r-1}x_0=x_{k_0}$ for some $k_0\in\N$, which contradicts the maximality of $x_0$. If $\frac t2<2^{r-1}x_0\le \frac {t2^{r-1}}{2^r-1}$, then
$$2^{r-1}x_0\pmod t=2^{r-1}x_0-t=x_{k_0}.$$
Consequently
$$x_{k_0}\le \frac {2^{r-1}t}{2^r-1}-t=-\frac{(2^{r-1}-1)t}{2^r-1}<0,$$
which yields a contradiction.

When $a=-2^{r-1}$, similarly we have $-2^{r-1}x_0+t:=x_{k_0}$. Then
$$t=2^{r-1}x_0+x_1\le \frac {t2^{r-1}}{2^r-1}+\frac t{2^r-1}=\frac{2^{r-1}+1}{2^r-1}t<t,$$
contradiction.

We finish the proof. \ez

\begin{lem}
	Let  $p\ne 2^r-1$ be an odd prime. If $O_{2}(p)$ is even and $r$ is odd,
	then $p^n$ is a complete number of the spectral pair $(\mu_{2^r, \{0, 2^{r-1}\}}, \Lambda(2^r, \{0, 1\}))$ for each $n\ge 1$.
\end{lem}
\proof   By definition  we have
$ 2^{rO_{2^r}(p^n)}\equiv1\pmod {p^n}$ for each $n\ge 1$. By Lemma \ref{lem6.1} and Proposition \ref{prop5.4}  we have
\begin{eqnarray}\label{eq6.1}
	O_{2^r}(p^n)=\frac{O_2(p^n)}{\gcd(O_2(p^n), r)}=\frac{O_2(p)p^{\max\{n-\ell_2(p), 0\}}}{\gcd(O_2(p)p^{\max\{n-\ell_2(p), 0\}}, r)},
\end{eqnarray}
which is even.
Write $rO_{2^r}(p^n)=\beta_n O_2(p)$, then $\beta_n$ is odd. By $ 2^{\beta_n O_2(p)}\equiv1\pmod {p^n}$ we have
\[(2^{\beta_n O_2(p)/2}-1)(2^{\beta_n O_2(p)/2}+1)\equiv0\pmod {p^n}.\]
If $2^{\beta_n O_2(p)/2}\equiv1\pmod p$, then $O_2(p)\mid \beta_n O_2(p)/2$, which is impossible. Then
$$(2^r)^{O_{2^r}(p^n)/2}=2^{\beta_n O_2(p)/2}\equiv-1\pmod {p^n}.$$
By Lemma \ref{lem6.2} the assertion follows. \ez
\begin{lem}
	Let  $p\ne 2^r-1$ be an odd  prime. If $O_{2}(p)$ is even and $r$ is even,
	then $p^n$ is a complete number of the spectral pair $(\mu_{2^r, \{0, 2^{r-1}\}}, \Lambda(2^r, \{0, 1\}))$ for each $n\ge 1$.
\end{lem}
\proof  Write $r=2^kr'$ and $O_2(p)=2^\al\al'$, where both $r'$ and $\al'$ are odd. By \eqref{eq6.1} we have
\begin{eqnarray*}
	rO_{2^r}(p^n)=\frac{2^{k+\al}r'\al'p^{\max\{n-\ell_2(p), 0\}}}{\gcd(2^\al \al'p^{\max\{n-\ell_2(p), 0\}}, 2^kr')}:=2^{k+\al-\min\{k, \al\}}\al'\beta_n,
\end{eqnarray*}
where $\beta_n$ is odd. When $\al\ge k$,  we have
\begin{eqnarray*}
	2^{2^{\al}\al'\beta_n}-1=2^{O_2(p)\beta_n}-1=(2^{O_2(p)\beta_n/2}-1)(2^{O_2(p)\beta_n/2}+1).
\end{eqnarray*}
Similar to the last proof, the assertion follows in this case.

When $\al< k$,  we have $ rO_{2^r}(p^n)=2^{k}\al'\beta_n=2^{k-\al}O_2(p)\beta_n$. Notice that
\begin{eqnarray*}
	2^{2^{k-\al}O_2(p)\beta_n}-1=\left( 2^{O_2(p)\beta_n}-1\right)\prod_{i=0}^{k-\al-1}\left( 2^{2^{i}O_2(p)\beta_n}+1\right)
\end{eqnarray*}
By the definition of $O_2(p)$ we obtain that $ 2^{O_2(p)\beta_n}\equiv1\pmod{p^n}$. Then with the same ideas we finish the proof. \ez
\begin{lem}
	Let  $p\ne 2^r-1$ be an odd prime and  $r\in 2\N$. If $O_{2}(p)$ is odd,
	then $p^n$ is a complete number of the spectral pair $(\mu_{2^r, \{0, 2^{r-1}\}}, \Lambda(2^r, \{0, 1\}))$ for each $n\ge 1$.
\end{lem}
\proof  Write $r=2^kr'$, where $k\ge 1$ and $r'$ is odd. By \eqref{eq6.1} we have
\begin{eqnarray*}
	rO_{2^r}(p^n)=\frac{2^{k}r'O_{2}(p)p^{\max\{n-\ell_2(p), 0\}}}{\gcd(O_{2}(p)p^{\max\{n-\ell_2(p), 0\}}, 2^kr')}:=2^kr'\beta_n,
\end{eqnarray*}
where $\beta_n$ is odd. According to $2^{2^kr'\beta_n}\equiv1\pmod{p^n}$, we have $4^{2^{k-1}r'\beta_n}\equiv1\pmod{p^n}$. Then
\[4^{2^{k-1}r'(\beta_n+1)}\equiv4^{2^{k-1}r'}\pmod{p^n}.\]
This implies that
\[\left(4^{2^{k-1}r'\frac{(\beta_n+1)}2}-2^{2^{k-1}r'}\right)\left(4^{2^{k-1}r'\frac{(\beta_n+1)}2}+2^{2^{k-1}r'}\right)\equiv0\pmod{p^n},\]
which is
$$(2^r)^{\frac{(\beta_n+1)}2}\equiv\pm 2^{2^{k-1}r'}\pmod{p^n}.$$
Note that $2^{2^{k-1}r'}<2^{r-1}-1$ for $r\ge 3$ and $2^{2^{k-1}r'}=2$ for $r=2$, then by Lemma \ref{lem6.2} and all $p^n$ are complete numbers of $(\mu_{2^r, \{0, 2^{r-1}\}}, \La(2^r, \{0, 1\})$ for all $n\ge 1$.
\ez

\medskip
The remaining case is that $rO_2(p)$ is odd, we guess that the same result is true but we cannot prove it. We can prove the following:
\begin{lem}\label{lem6.9}
	Let $p\neq 2^r-1$ be an odd prime. If $rO_2(p)$ is odd and  $r\nmid O_2(p)$, then  $p^n$ is complete number of the spectral pair $(\mu_{2^r, \{0, 2^{r-1}\}}, \Lambda(2^r, \{0, 1\}))$.
\end{lem}
To prove this lemma, we need a following proposition:
\begin{prop}\label{prop6.4}
	If $r=p^ks$ with $gcd(s, p)=1$, then $\ell_{2^r}(p)=\ell_2(p)+k$.
\end{prop}
\begin{proof}
	Fermat's little theorem indicates that $\gcd(p,O_2(p))=1$. By Lemma \ref{lem6.1}, we have
	\begin{eqnarray}\label{l2r1}
		O_{2^r}(p^{\ell_{2}(p)+k})&=&\frac{O_{2}(p^{\ell_{2}(p)+k})}{\gcd(O_{2}(p^{\ell_{2}(p)+k}),p^ks)}\nonumber\\
		&=&\frac{p^kO_{2}(p)}{\gcd(p^kO_{2}(p),p^ks)}=\frac{O_{2}(p)}{\gcd(O_{2}(p),r)}=O_{2^r}(p).
	\end{eqnarray}
	Similarly, we can obtain that
	\begin{equation}\label{l2r2}
		O_{2^r}(p^{\ell_{2}(p)+k+1})=\frac{p^{k+1}O_{2}(p)}{p^k\gcd(pO_{2}(p),s)}=\frac{pO_{2}(p)}{\gcd(O_{2}(p),r)}=pO_{2^r}(p).
	\end{equation}
	Then the conclusion follows from \eqref{l2r1}, \eqref{l2r2} and Proposition \ref{prop5.4}.
\end{proof}

{\noindent\bf Proof of Lemma \ref{lem6.9}.} We write $r=p^ks$ for some $k,s\in\N$ and $p\nmid s$. Then we have $\ell_{2^r}(p)=\ell_2(p)+k$ by Lemma \ref{prop6.4}.

Case 1: $k\ge 1$. For an odd prime $p$, we always have that
\begin{eqnarray}\label{eq6.5}
	2^{\frac{p-1}{2}}-1<\frac{1}{p}\left( 2^p-1\right)\le\frac{1}{p^k}\left( {2^{p^k}}-1\right)\le \frac{1}{p^k}\left( {2^{p^ks}}-1\right).
\end{eqnarray}
Fermat's little theorem shows that $O_2(p)\mid (p-1)$. Noticing that $O_2(p)$ and $p$ are odd numbers, then $O_2(p)\le\frac{p-1}{2}$.
By the definition of $\ell_{2}(p)$, we conclude that $p^{\ell_{2}(p)}\mid(2^{O_2(p)}-1)$. According to \eqref{eq6.5},
$$p^{\ell_{2^r}(p)}= p^k\cdot p^{\ell_{2}(p)}\le p^k\cdot \left( 2^{O_2(p)}-1\right)\le p^k\cdot \left( 2^{\frac{p-1}{2}}-1\right)< p^k\cdot \frac{1}{p^k}\left( {2^{p^ks}}-1\right)=2^r-1.$$
Corollary \ref{coro2.10} implies that $p^{\ell_{2^r}(p)}$ is a complete number. Then the result follows from Corollary \ref{coro5.8}.

Case 2: $k=0$ and $r\nmid O_2(p)$. We have
$$O_2(p^{\ell_{2^r}(p)})=O_2(p^{\ell_{2}(p)})=O_2(p)=\beta r+\alpha$$
for some $\alpha,\beta\in\N$ and $1\le\alpha\le r-1$. Then $2^{r(\beta+1)}\equiv 2^{r-\alpha}\pmod {p^{\ell_{2^r}(p)}}$ and $p^{\ell_{2^r}(p)}$ is a complete number. According to Corollary \ref{coro5.8}, We reach the conclusion.
\qed

In one word, we have proved
\begin{theo}
	Let $r\ge 2$ be an integer and let $p\ne 2^r-1$ be an odd prime. If $rO_2(p)$ is even or $r\nmid O_2(p)$, then $p^n$ is a complete number of the spectral pair $(\mu_{2^r, \{0, 2^{r-1}\}}, \Lambda(2^r, \{0, 1\}))$ for each $n\ge 1$.
\end{theo}

\end{document}